\documentclass{article}
\usepackage[english]{babel}
\usepackage[cp1251]{inputenc}
\usepackage{amsmath}
\usepackage{amssymb}
\usepackage{amsfonts}
\usepackage{amsthm}
\usepackage{graphicx}
\usepackage{url}

\def\thtext#1{
  \catcode`@=11
  \gdef\@thmcountersep{. #1}
  \catcode`@=12
}

\def\threst{
  \catcode`@=11
  \gdef\@thmcountersep{.}
  \catcode`@=12
}

 \catcode`@=11
 \def\.{.\spacefactor\@m}
 \catcode`@=12

\theoremstyle{plain}
\newtheorem{thm}{Theorem}
\newtheorem{cor}{Corollary}

\theoremstyle{definition}
\newtheorem{dfn}{Definition}

\newcommand{\cD}{\mathcal{D}}

\newcommand{\cM}{\mathcal{M}}

\newcommand{\cR}{\mathcal{R}}

\newcommand{\N}{\mathbb{N}}
\newcommand{\R}{\mathbb{R}}

\newcommand{\e}{\varepsilon}

\renewcommand{\r}{\rho}
\newcommand{\s}{\sigma}
\renewcommand{\v}{\varphi}

\newcommand{\dis}{\operatorname{dis}}

\renewcommand{\:}{\colon}

\renewcommand{\sp}{\supset}
\renewcommand{\ss}{\subset}
\newcommand{\x}{\times}
\def\rom#1{\emph{#1}}
\def\({\rom(}
\def\){\rom)}

\begin{document}
\title{Who Invented the Gromov--Hausdorff Distance?}
\author{Alexey A. Tuzhilin}
\maketitle

\begin{abstract}
One of the most beautiful notions of metric geometry is the Gromov--Hausdorff distance which measures the difference between two metric spaces. To define the distance, let us isometrically embed these spaces into various metric spaces and measure the Hausdorff distance between their images. The best matching corresponds to the least Hausdorff distance. The idea to compare metric spaces in such a way was described in M.~Gro\-mov publications~\cite{Gromov1}, \cite{Gromov2} dating back to 1981. It was shown that this distance being restricted to isometry classes of compact metric spaces forms a metric which is now called \emph{Gromov--Hausdorff metric}. However, whether M.~Gromov was the first who introduced this metric? It turns out that 6 years before these Gromov's works, in 1975, another mathematician, namely, David Edwards published a paper~\cite{Edwards} in which he defined this metric in another way. Also, Edwards found and proved the basic properties of the distance. It seems unfair that almost no one mentions this Edwards's paper including well-known specialists in metric geometry, see for instance~\cite{Gromov3}, \cite{Gromov4}, \cite{BurBurIvaEng}. The main goal of the present publication is to fill this gap.
\end{abstract}

In the present paper we talk about the authorship of the famous notion of metric geometry, namely, of the Gromov--Hausdorff distance. It is believed that this distance first appeared in the works of M.~Gromov~\cite{Gromov1}, \cite{Gromov2} dating back to 1981. However, in 1975 it was published a paper of D.~Edwards ``The Structure of Superspace'' \cite{Edwards} in which this Gromov--Hausdorff distance was defined in slightly different but equivalent terms.

The plan of this article is as follows. We start (sections~\ref{sec:GHdefProp}, \ref{sec:GHAlter}, and~\ref{sec:epsilonIsometry}) with well-known basic definitions and results related to the Gromov--Hausdorff distance, see~\cite{Gromov4}, \cite{BurBurIvaEng}, \cite{IvaNikolaevaTuz} for further details. Then we discuss (section~\ref{sec:Edwards}) which notions and results described above, and in which form, appeared in the pioneering work~\cite{Edwards} of D.~Edwards.

\section{Definition and Basic Properties of \\ the Gromov--Hausdorff Distance}\label{sec:GHdefProp}
\markright{\thesection.~Definition and Properties\\ of the Gromov--Hausdorff Distance}

Let $X$ be an arbitrary set. Each symmetric function $d:X\x X\to[0,\infty]$ we call a \emph{distance on $X$}. If $d$ is not equal to $\infty$, if it also is positively defined and satisfies the triangle inequality, then $d$ is called a \emph{metric}; a set $X$ with a metric $d$ defined on this $X$ is called a \emph{metric space}. If there is a distance function on $X$, then its values on a pair of points $x$, $y$ we denote by $|xy|$, if it does not cause confusion. If $A$ is a nonempty subset of $X$, and $x\in X$, then we put $|xA|=\inf\bigl\{|xa|:a\in A\bigr\}$. For $\e>0$ a subset $A$ of a metric space $X$ is called an \emph{$\e$-net\/} if $|xA|<\e$ for every $x\in X$.

A mapping $f\:X\to Y$ from a metric space $X$ to a metric space $Y$ is called \emph{isometric} if for any $x,x'\in X$ it holds $\bigl|f(x)f(x')\bigr|=|xx'|$. A bijective isometric mapping is called an \emph{isometry}. Denote by $\cM$ the set of all compact metric spaces considered upto an isometry.

Let $X$ be a metric space, $A$ a nonempty subset of $X$, and $r$ a positive real number. We define an \emph{$r$-neighborhood of $A$} as $U_r(A)=\{y\in X:|yA|<r\}$. If $A$ and $B$ are nonempty subsets of $X$, then the \emph{Gromov--Hausdorff distance between $A$ and $B$} is the value
$$
d_H(A,B)=\inf\bigl\{r:A\ss U_r(B)\ \&\ U_r(A)\sp B\bigr\}.
$$

Let $X$ and $Y$ be metric spaces. A triple $(X',Y',Z)$ consisting of a metric space $Z$ and its subsets $X'$ and $Y'$ isometric to $X$ and $Y$, respectively, is called a \emph{realization of the pair $(X,Y)$}. \emph{The Gromov--Hausdorff distance $d_{GH}(X,Y)$ between $X$ and $Y$} is the infimum of the numbers $r$ for which there exist realizations $(X',Y',Z)$ of the pair $(X,Y)$ with $d_H(X',Y')\le r$.

\begin{thm}\label{thm:GHdistanceIsMetric}
The distance function $d_{GH}$ restricted to $\cM$ is a metric.
\end{thm}

The set $\cM$ with the metric $d_{GH}$ is called the \emph{Gromov--Hausdorff space}. Now we present an initial segment of the properties list of the space $\cM$, which one usually meets, with some variations, in introduction sections to the subjects related to the Gromov--Hausdorff distance. Recall that a metric space $X$ is called \emph{geodesic}, if any points of  $X$ can be connected by a shortest curve which length equals to the distance between the points.

\begin{thm}\label{thm:GHproperties}
The Gromov--Hausdorff space
\begin{enumerate}
\item\label{thm:GHproperties:1} is contractible\/\rom;
\item\label{thm:GHproperties:2} is nonbounded and, thus, noncompact\/\rom;
\item\label{thm:GHproperties:3} contains the everywhere dense subset consisting of all finite metric spaces\/\rom;
\item\label{thm:GHproperties:4} is separable\/\rom;
\item\label{thm:GHproperties:5} is not locally compact\/\rom;
\item\label{thm:GHproperties:6} is complete\/\rom;
\item\label{thm:GHproperties:7} is geodesic.
\end{enumerate}
\end{thm}

\section{Other Definitions of \\ the Gromov--Hausdorff Distance}\label{sec:GHAlter}
\markright{\thesection.~Other Definitions of the Gromov--Hausdorff Distance}
The results we present below enable to give equivalent definitions of the Gromov--Hausdorff distance, which can be obtained from the right-hand sides of the formula in Theorem~\ref{thm:GH_admiss_metrics}, Theorem~\ref{thm:GH-distanceCorrespondences}, and Corollary~\ref{cor:GHdistanceIn TermsOfEllInfty}.

\subsection{Admissible Metrics on Disjoint Union of Spaces}
The next result means that to define the Gromov--Hausdorff distance it suffices to consider only metric spaces of the form $(X\sqcup Y,\r)$, where the restrictions of the metric $\r$ onto $X$ and $Y$ coincide with the corresponding initial metrics. We call such $\r$ an \emph{admissible metrics for $X$ and $Y$}, and the set of all admissible metrics for $X$ and $Y$ we denote by $\cD(X,Y)$. For each metric $\r\in\cD(X,Y)$ we denote by $\r_H$ the corresponding Hausdorff distance.

\begin{thm}\label{thm:GH_admiss_metrics}
For any metric spaces $X$ and $Y$ we have
$$
d_{GH}(X,Y)=\inf\bigl\{\r_H(X,Y):\r\in\cD(X,Y)\bigr\}.
$$
\end{thm}

\subsection{Correspondences and their Distortions}
Recall that for any sets $X$ and $Y$ each subsets of the Cartesian product $X\x Y$ is called a \emph{relation between $X$ and $Y$}. Let $\pi_X\:(x,y)\mapsto x$ and $\pi_Y\:(x,y)\mapsto y$ be the canonical projections. A relation $R$ is called a \emph{correspondence} if the restrictions of the projections $\pi_X$ and $\pi_Y$ onto $R$ are surjective. The set of all correspondences between $X$ and $Y$ is denoted by $\cR(X,Y)$.

If $X$ and $Y$ are some metric spaces, then for each nonempty relation $\s$ between $X$ and $Y$ its \emph{distortion\/} is
$$
\dis\s=\sup\Bigl\{\bigl||xx'|-|yy'|\bigr|:(x,y),(x',y')\in\s\Bigr\}.
$$

\begin{thm}\label{thm:GH-distanceCorrespondences}
For any metric spaces $X$ and $Y$ it holds
$$
d_{GH}(X,Y)=\frac12\inf\bigl\{\dis R:R\in\cR(X,Y)\bigr\}.
$$
\end{thm}

\subsection{Isometric Embeddings into the Space of\\ Bounded Sequences}
Denote by $\ell^\infty$ the space of all bounded real-values sequences with the metric generated by the norm $\bigl\|(x_1,x_2,\ldots)\bigr\|_\infty=\sup_{i\in\N}|x_i|$.

\begin{thm}\label{thm:Frechet-separable}
Let $X$ be a separable metric space, then there exists an isometric embedding of $X$ into $\ell^\infty$.
\end{thm}

\begin{cor}\label{cor:GHdistanceIn TermsOfEllInfty}
Let $X$ and $Y$ be separable metric spaces. Then
$$
d_{GH}(X,Y)=\inf d_H\bigl(\v(X),\psi(Y)\bigr),
$$
where the infimum is taken over all isometric embeddings $\v\:X\to\ell^\infty$ and $\psi\:Y\to\ell^\infty$.
\end{cor}

\section{The Gromov--Hausdorff Distance and\\ $\e$-Isometries}\label{sec:epsilonIsometry}
\markright{\thesection.~The Gromov--Hausdorff Distance and $\e$-isometries}

\begin{dfn}\label{dfn:epsilonIsometry}
For $\e>0$ a mapping $f\:X\to Y$ from a metric space $X$ to a metric space $Y$ is called an \emph{$\e$-isometry} if $\dis f\le\e$ and $f(X)$ is an $\e$-net.
\end{dfn}

\begin{thm}\label{thm:e-isometries}
Let $X$ and $Y$ be arbitrary metric spaces and $\e>0$. Then
\begin{enumerate}
\item\label{thm:e-isometries:1} if $d_{GH}(X,Y)<\e$, then there exists a $2\e$-isometry $f\:X\to Y$\rom;
\item\label{thm:e-isometries:2} if there exists an $\e$-isometry $f\:X\to Y$, then $d_{GH}(X,Y)<2\e$.
\end{enumerate}
\end{thm}

For compact metric spaces $X$ and $Y$ we define the distance $\widehat{d_{GH}}(X,Y)$ as the infimum of those $\e>0$ for which there exist $\e$-isometries $f\:X\to Y$ and $g\:Y\to X$. It is well-known that $\widehat{d_{GH}}(X,Y)$ does not coincide with $d_{GH}$.

\section{D.~Edwards Results}\label{sec:Edwards}
\markright{\thesection.~D.~Edwards Results}
Let us pass to description of the paper~\cite{Edwards}. D.~Edwards defined two distances between compact metric spaces. In the definition of the first distance he uses the notion of $\e$-isometry $f\:X\to Y$ which differs from the modern one presented in Definition~\ref{dfn:epsilonIsometry}.

\begin{dfn}
For metric spaces $X$, $Y$ and a real number $\e>0$ we call an \emph{$\e$-isometry in Edwards sense\/} each  mapping $f\:X\to Y$ with $\dis f\le\e$.
\end{dfn}

Thus, Edwards does not demand that the image $f(X)$ of $\e$-isometry $f\:X\to Y$ was everywhere dense in $Y$.

Further, Edwards defines a distance between compact metric spaces in the same manner like we did in section~\ref{sec:epsilonIsometry} for $\widehat{d_{GH}}$. To do that, he uses the notion of $\e$-isometry in his own sense. Notice that this distance plays an auxiliary role in Edwards constructions. Let us give a formal definition of this distance.

\begin{dfn}\label{dfn:EdwardsDistance1}
\emph{The first Edward's distance $d_E(X,Y)$} between compact metric spaces $X$ and $Y$ is the infimum of those $\e>0$ for which there exist $\e$-isometries $f\:X\to Y$ and $g\:Y\to X$ in the sense of Edwards.
\end{dfn}

Then Edwards shows that
\begin{enumerate}
\item the distance $d_E$ is a metric on the family $\cM$ of isometry classes of compact metric spaces;
\item the space $(\cM,d_E)$ is contractible;
\item finite metric spaces are everywhere dense in $(\cM,d_E)$;
\item the space $(\cM,d_E)$ is separable;
\item compact connected metric spaces are everywhere dense in the closed subspace of $(\cM,d_E)$ consisting of connected metric spaces;
\item there exists an isometric embedding of the ray $[0,\infty)\ss\R$ into the space $(\cM,d_E)$, thus the space $(\cM,d_E)$ is noncompact;
\item each Hilbert cube can be topologically embedded into $(\cM,d_E)$, thus the space $(\cM,d_E)$ is infinitely-dimensional;
\item the space $(\cM,d_E)$ does not have totally bounded neighborhoods, in particular, $(\cM,d_E)$ is not locally compact.
\end{enumerate}

Also, Edwards defines a second distance on the base of Theorem~\ref{thm:Frechet-separable}.

\begin{dfn}
\emph{The second Edward's distance $d_E^H(X,Y)$} between compact metric spaces  $X$ and $Y$ is the infimum of the Hausdorff distances between the sets $A,B\ss\ell^\infty$ isometric to $X$ and $Y$, respectively.
\end{dfn}

Notice that, by Corollary~\ref{cor:GHdistanceIn TermsOfEllInfty}, \textbf{the second Edward's distance coincides with the Gromov--Hausdorff distance\/} on the set $\cM$ of isometry classes of compact metric spaces.

Further, Edwards describes the properties of his second distance $d_E^H$:
\begin{enumerate}
\item the distance $d_E^H$ is a metric on the set $\cM$ of isometry classes of compact metric spaces (compare with Theorem~\ref{thm:GHdistanceIsMetric});
\item $d_E^H\ge\frac12d_E$ (similar with Theorem~\ref{thm:GH-distanceCorrespondences});
\item the space $(\cM,d_E^H)$ is contractible (compare with item~(\ref{thm:GHproperties:1}) of Theorem~\ref{thm:GHproperties});
\item finite metric spaces are everywhere dense in $(\cM,d_E^H)$ (compare with item~(\ref{thm:GHproperties:3}) of Theorem~\ref{thm:GHproperties});
\item the space $(\cM,d_E^H)$ is separable (compare with item~(\ref{thm:GHproperties:4}) of Theorem~\ref{thm:GHproperties});
\item compact connected polyhedra are everywhere dense in the closed subspace of  $(\cM,d_E^H)$ consisting of connected metric spaces;
\item there exists an isometric embedding of the ray $[0,\infty)\ss\R$ into the space $(\cM,d_E^H)$, thus the space $(\cM,d_E^H)$ is noncompact (compare with item~(\ref{thm:GHproperties:1}) of Theorem~\ref{thm:GHproperties});
\item the space $(\cM,d_E^H)$ does not have totally bounded neighborhoods, in particular, $(\cM,d_E^H)$ is not locally compact (compare with item~(\ref{thm:GHproperties:5}) of Theorem~\ref{thm:GHproperties});
\item the space $(\cM,d_E^H)$ is complete (compare with item~(\ref{thm:GHproperties:6}) of Theorem~\ref{thm:GHproperties}).
\end{enumerate}

Thus, Edwards in his paper~\cite{Edwards} \textbf{has defined the Gromov--Hausdorff distance\/} between isometry classes of compact metric spaces, and \textbf{laid down the foundations of the relevant theory by proving a number of results describing the properties of this distance}. Taking into account the above discussion, we see that the work~\cite{Edwards} of David Edwards is a classical one in metric geometry field, thus, the absence of references to this work is insulting undeserved accident. I think that this work have to take its rightful place of honor in the history of metric geometry, and I hope that the present paper, both with my corrections of Russian and English Wikipedia, will support this idea.

\section{Addendum: the known citations of the D.~Edwards paper}
\markright{\thesection.~Addendum: the known citations of the D.~Edwards paper}
To complete the present paper, let us present the citations of the paper~\cite{Edwards} which I found in literature (of course, I do not pretend to be complete).

\textbf{Wikipedia\/} (before my correction): Russian Wikipedia did not have any mentioning of the Edwards's work. In English Wikipedia is was written the following: ``The notion of Gromov--Hausdorff convergence was first used by Gromov to prove that any discrete group with polynomial growth is virtually nilpotent (i.e. it contains a nilpotent subgroup of finite index). See Gromov's theorem on groups of polynomial growth. (Also see D. Edwards for an earlier work.)''

The only author wrote about topics related to the Gromov--Hausdorff distance and mentioned the Edward's paper, as I found, was \textbf{F.Latremoliere}, see for example~\cite{Latremoliere}. Here it was written the following: ``On the other hand, Edward~\cite{Edwards} introduced what is now known as the Gromov--Hausdorff distance for compact metric spaces, based upon consideration about quantum gravity, suggesting that one would have to work, not on one spacetime, but rather on a “super-space” of possible space--times, metrized by a generalized Hausdorff distance, thus allowing fluctuations in the gravitational field. Edwards was following upon the suggestion to study quantum gravity by means of a superspace by Wheeler~\cite{Wheeler}. It should be noted that Gromov introduces his metric~\cite{Gromov2} for more general locally compact metric spaces, and used his metric toward major advances in group theory.''

On the \textbf{David Edwards cite}~\cite{EdwardsCite}, the author of~\cite{Edwards} has written: ``This paper is related to Wheeler's notions of quantum geometrodynamics; it includes my defining a complete metric on the space of isometry classes of compact metric spaces --- it was later rediscovered by Gromov and is called the Gromov metric.''

\vfill\eject

\end{document}